\documentclass[11pt]{article}
\usepackage{amssymb}
\usepackage{soul}
\usepackage{color, xcolor}

\usepackage{latexsym}
\usepackage{amsfonts}
\def\no{\nonumber}
\def\be{\begin{equation}}
\def\ee{\end{equation}}
\def\ba{\begin{eqnarray}}
\def\ea{\end{eqnarray}}

\def\btd{\bigtriangledown}
\def\e1{\epsilon}
\def\AAl{\mathcal{A}_{\lambda}}
\def\A0{\stackrel{\circ}{\AAl}}

\def\o1{\omega}
\def\01{\Omega}
\def\c1{\gamma}
\def\a{\alpha}
\def\g1{\Sigma}
\def\bigcap{\cap}
\def\bigcup{\cup}

\def\l1{\Lambda}

\def\v1{\varphi}

\def\d1{\delta}
\def\part{\partial}

\def\f1{\frac}
\def\t1{\theta}
\def\b1{\beta}

\def\bs{\begin{eqnarray*}}
\def\es{\end{eqnarray*}}
\def\bes{\begin{equation*}}
\def\ees{\end{equation*}}
\def\m1{\Theta}
\def\w1{\wedge}

\date{}

\begin{document}
\large{\title {A blow-up formula for stationary quaternionic maps
\thanks{This work is supported by NSF grant 11721101.}
}
\author{Jiayu Li\thanks{{jiayuli@ustc.edu.cn}}\\
\and
Chaona Zhu\thanks{zcn1991@mail.ustc.edu.cn}}

\maketitle

\newtheorem{theorem}{Theorem}[section]
\newtheorem{lemma}[theorem]{Lemma}
\newtheorem{corollary}[theorem]{Corollary}
\newtheorem{remark}[theorem]{Remark}
\newtheorem{definition}[theorem]{Definition}
\newtheorem{proposition}[theorem]{Proposition}


\begin{abstract}
Let $(M, J^\alpha, \alpha =1,2,3)$ and $(N, {\cal J}^\alpha,
\alpha =1,2,3)$ be Hyperk\"ahler manifolds. Suppose that $u_k$ is a sequence of
stationary quaternionic maps and converges weakly to $u$ in $H^{1,2}(M,N)$, we derive a blow-up formula for
$\lim_{k\to\infty}d(u_k^*{\cal J}^\alpha)$, for $\alpha=1,2,3$, in the weak sense. As a corollary, we show that the maps constructed by
Chen-Li [CL2] and by Foscolo [F] can not be tangent maps (c.f [LT], Theorem 3.1) of a stationary quaternionic map satisfing $d(u^*{\cal J}^\alpha)=0$.

\end{abstract}

\section{Introduction and the main result}
A hyperk\"ahler manifold is a Riemannian manifold $(M,g)$ with three parallel
complex structures $\left\{J^1,J^2,J^3\right\}$ compatible with the metric $g$ such that
$(J^1)^2=(J^2)^2=(J^3)^2=J^1J^2J^3=-id$.
The simplest hyperk\"ahler manifold is the Euclidean space $\mathbb{R}^{4m}$.
It is well-known that the only compact hyperk\"ahler manifolds of dimension $4$
are $K3$ surfaces and complex tori. Let $(M, g, J^\alpha, \alpha =1,2,3)$ and $(N, h, {\cal J}^\alpha,
\alpha =1,2,3)$ be hyperk\"ahler manifolds.  Let $\omega_\alpha(\cdot ,\cdot)=g(\cdot ,J^\alpha\cdot)$ and $\Omega_\alpha(\cdot ,\cdot)=h(\cdot ,{\cal J}^\alpha\cdot)$, $(\alpha =1,2,3)$ be the K\"ahler forms on M and N respectively.
A smooth map
$u:M\rightarrow N$ is called a quaternionic map (triholomorphic map) if
\begin{equation}\label{trihol}
A_{\a\beta}{\cal J}^{\beta}\circ du\circ J^{\a}=du
\end{equation}
where $A_{\a\beta}$ denote the entries of a matrix $A$ in $SO(3)$.
For simplicity, we choose $A_{\a\beta}=\delta_{\a\beta}$.

The quaternionic maps (triholomorphic maps) between Hyperk\"ahler manifolds has been studied
by many aothors (cf. [BT], [Ch], [CL1}, [CL2], [FKS], [W]).
Quaternionic maps automatically minimize the energy functional in
their homotopy classes (cf. [Ch], [CL1] and
[FKS]) and hence they are harmonic. It can be verified that holomorphic and anti-holomorphic
maps with respect to some complex structures on $M$ and $N$ are quaternionic maps. However, Chen-Li constructed
quaternionic maps which are not holomorphic with respect to any complex structures on $M$ and $N$ (cf. [CL1]).

\begin{definition}
A map
$u$ from $M$ to $N$ is called a {\it stationary quaternionic map} if
it is a stationary harmonic map and it is a quaternionic map outside
its singular set.
\end{definition}

It is clear that (c.f. [BT]), if $u$ satisfies (\ref{trihol}) almost everywhere, and
\begin{equation}\label{QMesid}
d(u^*{\cal J}^\alpha)=0,~~{\rm for}~~\alpha=1,2,3,
\end{equation}
then $u$ is a stationary quaternionic map.

Chen-Li ([CL2]) proved that, if there is a harmonic sphere $\phi: {\mathbb S}^2\to N$ which satisfies
\begin{equation}\label{eqn:1}
d\phi\,J_{{\mathbb S}^2} =-\sum^3_{k=1}a_k {\cal J}^k\, d\phi,
\end{equation}
where $\vec{a}=(a_1,a_2,a_3): {\mathbb S}^2\to {\mathbb S}^2$, and
\begin{equation}\label{stationary}
\int_{{\mathbb S}^2}x_i|\nabla \phi|^2d\sigma=0,\,\,\,i=1,2,3,\,\,\,
(x_1,x_2,x_3)\in {\mathbb S}^2,
\end{equation}
then
$$
u(x,x^4)=\phi(\frac{x}{|x|}) ~~{\rm for~~ any} ~x\in \mathbb{R}^3\backslash \{0\}
$$
is a stationary quaternionic map with the $x^4$-axis as its singular set.

Chen-Li ([CL2]) showed that there does exist a complete noncompact hyperk\"ahler manifold, into which there is a harmonic ${\mathbb S}^2$ which satisfies
(\ref{eqn:1}) and (\ref{stationary}). Recently, Foscolo [F] showed that there exists a compact $K3$ surface with the above property. However, the map $u$ constructed by Chen-Li or by Foscolo does not satisfy (\ref{QMesid}). Now the question is whether the maps constructed by Chen-Li or by  Foscolo could be a tangent map of a stationary quaternionic map with identity (\ref{QMesid}), if not the singular set of a stationary quaternionic map with identity (\ref{QMesid}) might be of codimensional 4 (Remark 1.2 in [BT]).

Suppose that $u_k$ is a sequence of
stationary quaternionic maps with bounded energies $E(u_k)\leq \Lambda$.
The {\it blow-up set of $u_k$} can be defined as
$$
\g1 =\bigcap_{r>0}\{x\in M |
\liminf _{k\to\infty}r^{2-m}\int_{B_r(x)}|\btd u_k|^2dy\geq \e1 _0\}.
$$
We can always assume that $u_k\rightharpoonup u$ weakly in $W^{1,2}(M,N)$
and that
$$
|\btd u_k|^2dx\rightharpoonup |\btd u|^2dx+\nu
$$
in the sense of measure as $k\to\infty$. Here $\nu$ is a nonnegative
Radon measure on $M$ with support in $\g1$. It is known that $\g1$ is
a ${\cal H}^{m-2}$-rectifiable set, and we may write $\nu =\theta (x)H^{m-2}\lfloor \Sigma$.
It is clear that strongly convergence in $H^{1,2}(M,N)$ preserves the identity (\ref{QMesid}). In this paper we mainly prove the following
blow-up formula for weakly convergence sequence of stationary quaternionic maps.

\begin{theorem}\label{iden}
  Let $u_k$ be a sequence of stationary quaternionic map with $E(u_k)\leq\Lambda$. Assume that $u_k\to u$ weakly in $H^1(M,N)$. Then there exist $\ (a^1, a^2, a^3)\in\mathbb{R}^3$ with $\sum_{\alpha=1}^3(a^\alpha)^2=1$ such that, for any smooth $(m-3)$-form $\eta$  with compact support in $M$,
  \be\label{stafor}
  \lim_{k\to\infty}\sum_{\alpha=1}^3a^\alpha\int_{M}d\eta\wedge u_k^*{\cal J}^\alpha=\sum_{\alpha=1}^3a^\alpha\int_{M}d\eta\wedge u^*{\cal J}^\alpha+\int_{\Sigma}\theta d\eta|_\Sigma
  \ee
  and for any $(b^1, b^2, b^3)\perp(a^1, a^2, a^3)$, there holds
  \bs
  \lim_{k\to\infty}\sum_{\alpha=1}^3b^\alpha\int_{M}d\eta\wedge u_k^*{\cal J}^\alpha=\sum_{\alpha=1}^3b^\alpha\int_{M}d\eta\wedge u^*{\cal J}^\alpha.
  \es
\end{theorem}

As a corollary of the theorem, the maps constructed by
Chen-Li [CL2] and by Foscolo [F] can not be tangent maps (c.f [LT], Theorem 3.1) of a stationary quaternionic map satisfing $d(u^*{\cal J}^\alpha)=0$.

\section{The proof of the blow-up formula}

If $u$ is a strong limit of a sequence
of stationary quaternionic maps in $H^{1,2}(M,N)$, then it's easy to see that $u$ satisfies (\ref{QMesid}).
If $u$ is just a weak limit, i.e. there exists a sequence of stationary quaternionic maps $u_k$ satisfying $u_k\to u$ weakly in $H^{1,2}(M, N)$ and $|\nabla u_k|^2dV\to|\nabla u|^2dV+\theta H^{m-2}|_{\Sigma}$ in the sense of measure, we prove in this section a formula for the blow-up set $\theta H^{m-2}|_{\Sigma}$ and the limiting map $u$.

Without loss of generality, we may assume that $m=4$. Because $\g1$ is a
$H^{m-2}$-rectifiable set, so we may assume that
$\g1 =\bigcup _{i=0}^{\infty}\g1 _i$, $\g1_i\cap\g1_{i'}
=\phi$ if $i\neq i'$, $H^{m-2}(\g1_0)=0$, $\g1_i
\subset N_i$ and $N_i$ $(i=1,2,\cdots )$ is an
$(m-2)$-dimensional embedded $C^1$ submanifold of $M$.
It is important that (see p. 61
in [Si])
$T_x\g1 = T_x N_i$ for $H^{m-2}$-a.e. $x\in\g1_i$.

It is known that $\nu =\theta (x)H^{m-2}\lfloor \g1$,
where $\theta (x)$ is upper semi-continuous with
$\e1_0\leq \theta (x)\leq C(\l1)$ for $H^{m-2}$-a.e.
$x\in \g1$, $C(\l1)$ is a positive constant depending only on $M$
and $\l1$ (cf. [Lin], Lemma 1.6).
Since
$H^{m-2}(\g1)<+\infty$, for any $\d1 >0$,
there exist
$\g1_{\d1}\subset \g1$ and  $i_0$
such that
$H^{m-2}(\g1_{\d1})<\d1$, $\g1^c_{\d1}=\g1\setminus \g1_{\d1}
= \cup _{i=1}^{i_0} \g1^{\d1}_i$
where $ \g1^{\d1}_i\subset \g1_i$ $(i=1,\cdots ,i_0)$ is a bounded
closed set.
We choose a covering $\{ B_{r_n} | n=1,2,\cdots\}$ of $\g1_{\d1}$
such that $\sum_{n}r_n^{m-2}< C\d1$.
Here and in the sequel, $C$ always denotes a uniform constant depending
only on $M$ and $N$.

Suppose that $(x^1, ...,x^{4})$ is a local normal coordinate system in $B_{\e1}(\Sigma_{i}^\delta)$, and that $(x^3,x^{4})$ is the corresponding coordinate system in $\Sigma_i$,
 and
the matrix expressions of the complex structures are given by
(\ref{J11}), (\ref{J22}) and (\ref{J33}).

\be\label{J11}
J^1=\left ( \begin{array}{cccc}
0&0&0&-1\\
0&0&1&0\\
0&-1&0&0\\
1&0&0&0
\end{array}\right ),\,\,\,\,
A_{1\beta}{\cal J}^\beta=\left ( \begin{array}{ccccc}
J^1& & &  \\
&\cdot & &  \\
& & \cdot  & \\
& & & J^1
\end{array}\right )
\ee
\be\label{J22}
J^2=\left ( \begin{array}{cccc}
0&-1&0&0\\
1&0&0&0\\
0&0&0&1\\
0&0&-1&0
\end{array}\right ),\,\,\,\,
A_{2\beta}{\cal J}^\beta=\left ( \begin{array}{ccccc}
J^2& & &  \\
&\cdot & &  \\
& & \cdot & \\
& & & J^2
\end{array}\right )
\ee
\be\label{J33}
J^3=\left ( \begin{array}{cccc}
0&0&1&0\\
0&0&0&1\\
-1&0&0&0\\
0&-1&0&0
\end{array}\right ),\,\,\,\,
A_{3\beta}{\cal J}^\beta=\left ( \begin{array}{ccccc}
J^3& & &  \\
&\cdot &  & \\
& & \cdot  & \\
& & &  J^3
\end{array}\right )
\ee
where $A_{\a\beta}{\cal J}^\beta$ are $4n\times 4n$-matrices, $A_{\a\beta}$ are the entries of a matrix $A$ in $SO(3)$.
Then the quaternionic equation is
\begin{equation}\label{qmap11}
\left\{\begin{array}{clcr}
u^1_1+u^2_2+u^3_3+u^4_4&=&0\\
u^2_1-u^1_2+u^4_3-u^3_4&=&0\\
u^3_1-u^1_3-u^4_2+u^2_4&=&0\\
u^4_1-u^1_4-u^2_3+u^3_2&=&0\\
u^5_1+u^6_2+u^7_3+u^8_4&=&0\\
u^6_1-u^5_2+u^8_3-u^7_4&=&0\\
u^7_1-u^5_3-u^8_2+u^6_4&=&0\\
u^8_1-u^5_4-u^6_3+u^7_2&=&0\\
\cdots\\
u^{4n-3}_1+u^{4n-2}_2+u^{4n-1}_3+u^{4n}_4&=&0\\
u^{4n-2}_1-u^{4n-3}_2+u^{4n}_3-u^{4n-1}_4&=&0\\
u^{4n-1}_1-u^{4n-3}_3-u^{4n}_2+u^{4n-2}_4&=&0\\
u^{4n}_1-u^{4n-3}_4-u^{4n-2}_3+u^{4n-1}_2&=&0.
\end{array}\right.
\end{equation}

\begin{theorem}\label{iden-simple}
 For any smooth $(m-3)$-form $\eta$ with compact support in $M$, we have
  \bs
  \lim_{k\to\infty}\sum_{\alpha=1}^3A_{\a\beta}\int_{M}d\eta\wedge u_k^*{\cal J}^\beta=\sum_{\alpha=1}^3A_{\a\beta}\int_{M}d\eta\wedge u^*{\cal J}^\beta+\int_{\Sigma}\theta d\eta|_\Sigma
  \es
  and
  \bs
  \lim_{k\to\infty}A_{1\beta}\int_{M}d\eta\wedge u_k^*{\cal J}^\beta=A_{1\beta}\int_{M}d\eta\wedge u^*{\cal J}^\beta,
  \es
  \bs
  \lim_{k\to\infty}A_{3\beta}\int_{M}d\eta\wedge u_k^*{\cal J}^\beta=A_{3\beta}\int_{M}d\eta\wedge u^*{\cal J}^\beta,
  \es
\end{theorem}
{\it Proof.} Assume that $\eta=\sum_{I}\eta_Idx^I$. We have
\begin{eqnarray}
\lefteqn{\lim_{k\to\infty}\int_{M}d\eta\wedge u_k^*(A_{\a\beta}{\cal J^\beta})=\int_{M}d\eta\wedge u^*( A_{\a\beta}{\cal J}^\beta)}~~~~~~~~~~~~~~~~~~~~~~~~~\no\\
&+&\lim_{\delta\to0}\lim_{\epsilon\to 0}\lim_{k\to\infty}\int_{B_{\epsilon}(\bigcup_{i=1}^{i_0}\Sigma_i^\delta)}d\eta\wedge u_k^*( A_{\a\beta}{\cal J}^\beta)~~~~~~~~~~~~~~~~~~~~~~~~~\no\\
&+&\lim_{\delta\to0}\lim_{\epsilon\to0}\lim_{k\to\infty}\int_{\bigcup_{n}B_{r_n}\setminus B_{\epsilon}(\bigcup_{i=1}^{i_0}\Sigma_i^\delta)}d\eta\wedge u_k^*( A_{\a\beta}{\cal J}^\beta).
\end{eqnarray}
It's easy to see that
\be
\lim_{\delta\to0}\lim_{\epsilon\to0}\lim_{k\to\infty}\int_{\bigcup_{n}B_{r_n}}d\eta\wedge u_k^*({\cal J}^\beta)=0
\ee

By Lemma 2.2 in [LT], we get
\ba\label{for1}
\lefteqn{\lim_{\delta\to0}\lim_{\epsilon\to0}\lim_{k\to\infty}\int_{B_{\epsilon}(\Sigma_i^\delta)}d\eta\wedge u_k^*(A_{\a\beta}{\cal J}^\beta)}\no\\
&=&\lim_{\delta\to0}\lim_{\epsilon\to0}\lim_{k\to\infty}\int_{B_{\epsilon}(\Sigma_i^\delta)}2\f1{\partial{\eta_{I}}}{\partial x^l}\f1{\partial u_k^\sigma}{\partial x^1}(A_{\a\beta}{\cal J}^\beta)_{\sigma\gamma}\f1{\partial u_k^\gamma}{\partial x^2}dx^l\wedge dx^{I}\wedge dx^1\wedge dx^2
\ea
Substituting (\ref{qmap11}) to (\ref{for1}) and applying Lemma 2.2 in [LT], we have
\bs
\lim_{\delta\to0}\lim_{\epsilon\to0}\lim_{k\to\infty}\int_{B_{\epsilon}(\Sigma_i^\delta)}d\eta\wedge u_k^*(A_{1\beta}{\cal J}^\beta)=\lim_{\delta\to0}\lim_{\epsilon\to0}\lim_{k\to\infty}\int_{B_{\epsilon}(\Sigma_i^\delta)}d\eta\wedge u_k^*(A_{3\beta}{\cal J}^\beta)=0
\es
and
\ba\label{for2}
\lefteqn{\lim_{\delta\to0}\lim_{\epsilon\to0}\lim_{k\to\infty}\int_{B_{\epsilon}(\Sigma_i^\delta)}d\eta\wedge u_k^*(A_{2\beta}{\cal J}^\beta)
=\lim_{\delta\to0}\lim_{\epsilon\to0}\lim_{k\to\infty}\int_{B_{\epsilon}(\Sigma_i^\delta)}|\nabla u_k|^2d\eta\wedge dx^1\wedge dx^2}\no\\
&=&\lim_{\delta\to0}\lim_{\epsilon\to0}(\int_{B_{\epsilon}(\Sigma_i^\delta)}|\nabla u|^2d\eta\wedge dx^1\wedge dx^2+\int_{B_{\epsilon}(\Sigma_i^\delta)\bigcap\Sigma}\theta d\eta|_{\Sigma})
=\int_{\Sigma_i}\theta d\eta|_{\Sigma}
\ea
Then the proof of the theorem is completed.
\hfill Q.E.D.

\begin{remark}
From this theorem, we see that if $u_k$ satisfies (\ref{QMesid}), the weak limit $u$ still satisfies (\ref{QMesid}) if and only if $\theta={\rm constant}$.
\end{remark}

As a corollary, we can derive that $\theta(x)$ is locally constant. Precisely,

\begin{corollary}
Under the assumption of Theorem \ref{iden}, and assume that there is an open ball $B^{m}\subset M\setminus{\rm Sing_u}$ with $H^{m-2}(\Sigma\cap B^m)>0$. We have $\theta(x)$ is constant on $\Sigma\cap B^m$.
\end{corollary}
{\it Proof.}
In (\ref{stafor}), we choose cutoff function $\eta$ such that ${\rm supp}\eta\subset B^{m}$. Since $B^m\subset M\setminus{\rm Sing_u}$, we have $u$ is smooth on $B^m$. Then $du^*{\cal J}^\beta=0$ on $B^m$ for $\beta=1,2,3$. In view of (\ref{stafor}), we conclude that $\theta$ is constant on $\Sigma\cap B^m$.

\hfill Q.E.D.

\

Let $\phi~:~{\mathbb S}^2\to N$ be a nonconstant smooth map satisfying (\ref{eqn:1}) and (\ref{stationary}).  Set
\begin{equation}\label{tc}
u(x,x^4)=\phi(\frac{x}{|x|}) ~~{\rm for~~ any} ~~x\in \mathbb{R}^3\backslash \{0\}~~x^4\in\mathbb{R}^{m-3}
\end{equation}
as Chen-Li ([CL2]) did.  
Then we have
\begin{proposition}\label{yn}

For any smooth $(m-3)$-form $\eta$ with compact support in $\mathbb{R}^{m}$, we have
\begin{equation}\label{yniden}
\int_{\mathbb{R}^{m}}d\eta \wedge u^*{\cal J}^\alpha=-E_T^\alpha (\phi)\int_{\mathbb{R}^{m-3}}\eta(0,x^4),
\end{equation}
where
$$
E_T(\phi)= \int_{{\mathbb S}^2} \langle J^\a_{{\mathbb S}^2},
u^*{\cal J}^{\a}\rangle d\sigma .
$$
\end{proposition}
{\it Proof.}
We choose a spherical coordinate system $(r, \varphi ,\theta )$ in $\mathbb{R}^3$, because $u$ is smooth for any $r>0$, we have
\bs
\lefteqn{\int_{\mathbb{R}^{m}}d\eta \wedge u^*{\cal J}^\alpha}\\
&=&\int_{\mathbb{R}^{m-3}} \int_0^\infty \frac{\partial \eta_I }{\partial r}dr\wedge dx^I\int_{{\mathbb S}^2}  \phi^*{\cal J}^\alpha\\
&=&-\int_{\mathbb{R}^{m-3}}  \eta (0,x^4)\int_{{\mathbb S}^2}  \phi^*{\cal J}^\alpha\\
&=&-E_T^\alpha (\phi)\int_{\mathbb{R}^{m-3}}\eta(0,x^4)
\es

\hfill Q.E.D.

By Theorem \ref{iden-simple} and Proposition \ref{yn}, we have the following corollary.

\begin{corollary}
The map $u$ defined in (\ref{tc}) can not be a tangent map (c.f [LT], Theorem 3.1) of a stationary quaternionic map with the property (\ref{QMesid})
at a singular point.
\end{corollary}

{\it Proof.}
Suppose that $u$ is defined as in (\ref{tc}). If it is a tangent map, then we have by Theorem \ref{iden-simple},
$$
\sum_{\alpha=1}^3A_{\a\beta}\int_{M}d\eta\wedge u^*{\cal J}^\beta+\int_{\Sigma}\theta d\eta|_\Sigma=0.
$$
By
Proposition \ref{yn}, we obtain
$$
\sum_{\alpha=1}^3A_{\a\beta}E_T^\beta (\phi)\int_{\mathbb{R}^{m-3}}\eta(0,x^4)=\int_{\Sigma}\theta d\eta|_\Sigma.
$$
Since $u$ is stationary, by the blow-up formula of Li-Tian [LT], we have $\Sigma$ is stationary. Using the constancy theorem (Theorem 41.1 in [Si]), it follows that the density function $\theta$ is constant in every connected component of $\Sigma$, which implies that $\phi$ is homotopy to a constant map. We therefore get a contradiction.

\hfill Q.E.D.

\begin{center}
{\large\bf REFERENCES}
\end{center}
\footnotesize
\begin{description}

\item[{[BT]}] {C. Bellettini and G. Tian, Compactness results for triholomorphic maps, J. Eur. Math. Soc., 2(2019), 1271-1317.}
\item[{[Ch]}] {J. Chen, Complex anti-self-dual connections on
product of Calabi-Yau surfaces and triholomorphic curves, Commun.
Math. Phys. 201(1999), 201-247.}
\item[{[CL1]}] {J. Chen and J. Li, Quaternionic maps between Hyperk\"ahler
manifolds, J. Diff. Geom. 55(2000), no. 2, 355-384.}
\item[{[CL2]}] {J. Chen and J. Li, Quarternionic maps and minimal surfaces, Ann. Sc. Norm. Super. Pisa Cl. Sci. 4 (2005), no. 3, 375-388.}
\item[{[FKS]}] {J.M. Figuroa-O'Farrill, C. K\"ohl and B. Spence,
Supersymmetric Yang-Mills, octonionic instantons and triholomorphic
curves, Nucl. Phys. B 521 (1998) no. 3, 419-443.}
\item[{[F]}] {L. Foscolo, ALF gravitational instantons and collapsing Ricci-flat metrics on the K3 surface, J. Diff. Geom., 112(2019), 79-120.}
\item[{[LT]}] {J. Li, and G. Tian,
A blow-up formula for stationary harmonic maps, IMRN, 14(1998), 735-755.}
\item[{[Lin]}] {F.-H. Lin, Gradient estimates and blow-up analysis
for stationary harmonic maps I, Ann. of Math. 149(1999),
785-829.}
\item[{[Si]}] {L. Simon,  Lectures on Geometric Measure Theory,
Proc. Center Math. Anal. 3(1983), Australian National Univ.
Press.}
\item[{[W]}] {C. Wang,  Energy quantization for triholomorphic maps, Calc. Var. PDE 18(2003), 145-158.}

\end{description}

\end{document}